\documentclass{amsart}
\input{diagrams.tex}
\diagramstyle[midshaft,scriptlabels,height=8mm,width=8mm]
\def\Sec{{\mathop{\rm Sec}\nolimits\,}}
\def\Tan{{\mathop{\rm Tan}\nolimits\,}}
\def\Hilb{{\mathop{\rm Hilb}\nolimits\,}}

\def\cN{{\mathcal N}}

\def\bbbp{{\mathbb P}}
\def\bbbc{{\mathbb C}}

\def\bbbg{{\mathbb G}}

\def\hto{\hookrightarrow}

\theoremstyle{definition}
\newtheorem{defn}{Definition}[section]

\theoremstyle{plain}

\newtheorem{theorem}[defn]{Theorem}
\newtheorem{lemma}[defn]{Lemma}
\newtheorem{corollary}[defn]{Corollary}

\theoremstyle{remark}
\newtheorem{remark}[defn]{Remark}

\begin{document}

\title{A Note on Generic Projections}

\author{Hubert Flenner}
\address{Fakult\"at f\"ur Mathematik der Ruhr-Universit\"at,
Universit\"atsstr.\ 150, Geb.\ NA 2/72, 44780 Bochum, Germany}
\email{Hubert.Flenner@ruhr-uni-bochum.de}

\author{Mirella Manaresi}
\address{Dipartimento di Matematica, Universit\`a\\
P.za di Porta S.\ Donato 5,
40126 Bologna, Italy}
\email{manaresi@dm.unibo.it}

\thanks{This paper was supported by  the DFG Schwerpunkt "Global
methods in complex analysis", MIUR and the Research Group GNSAGA of INDAM.
This investigation was also supported by the University of Bologna, funds for selected
research topics}
%\Names{Hubert Flenner, Mirella Manaresi}
\maketitle

\begin{abstract}
Let $X \subseteq \bbbp^N = \bbbp^{2n}_K$ be a subvariety of dimension $n$
and $P \in \bbbp^N$ a generic point. If the tangent variety $\Tan
X$ is equal to $\bbbp^N$ then for generic points $x$, $y$ of $X$ the
projective tangent spaces
$t_xX$ and $t_yX$ meet in one point $P=P(x,y)$. The main result of this paper is
that the rational map $(x,y)\mapsto P(x,y)$ is dominant. In other words,  a generic
point $P$ is uniquely determined by the ramification locus
$R(\pi_P)$ of the linear projection $\pi_P:X\to \bbbp^{N-1}$.
\end{abstract}

\section{Introduction}
Let $X \subseteq \bbbp^N = \bbbp^N_K$ be a subvariety of dimension $n$ and
$\Lambda \subseteq \bbbp^N$ a linear subspace of codimension $k+1$ with $k \ge
n$. In \cite{FM2} we studied the question how the ramification locus
$R(\pi_\Lambda)$ of the linear projection $\pi_\Lambda:X\to \bbbp^k$ with center
$\Lambda$ varies with $\Lambda$. More precisely we considered the problem as to how
$\Lambda$ is determined by $R(\pi_\Lambda)$.

The motivation for this question comes from the study of the St\"uckrad-Vogel cycle.
If we denote by $v=v(X,X)=\sum_iv_i$ the St\"uckrad-Vogel selfintersection cycle of
$X$ then by \cite{vGa, FM1} the cycle $v_{2n-k-1}$ may be interpreted as a
ramification cycle of the generic projection $\pi_{\Lambda_U}:X\to\bbbp^{k}$, where
$\Lambda_U\subseteq
\bbbp^N_L$ is the generic linear subspace of codimension $k+1$ given by the
equations $\sum_jU_{ij}x_j=0$, $0\le i\le k$, over the purely transcendental
extension field
$L:=K(U_{ij})$. 

In \cite{FM2} we proved that the cycle $v_{2n-k-1}$ has maximal transcendence
degree if and only if the rational map $\Lambda\mapsto R(\pi_\Lambda)$ from the
Grassmannian
$\bbbg(N-k-1,N)$ of $(N-k-1)$-planes in $\bbbp^N$ to the Hilbert scheme $\Hilb_X$ of
$X$ is generically finite. In particular we treated in that paper the case of smooth
surfaces in
$\bbbp^4$ and the projection from a point, i.e.\ $n=2$, $N=4$ and $\Lambda=\{P\}$.
Under a suitable positivity condition on the normal bundle
$\cN_{X/\bbbp^4}(-1)$ we showed that $R(\pi_P)$ uniquely determines $P$, and from
this we deduced that $v_0(X,X)$ has maximal transcendence degree. 

In a letter to the second author C.\ Ciliberto improved this by
showing that for any nondegenerate smooth surface $X$ in $\bbbp^4$ the same result
remains true. His essential idea was to use the second fundamental form of the
surface to conclude that the rational map
$$
\begin{diagram}[height=15pt]
X\times X&\rDashto&\bbbp^4\\
(x,y)& \rMapsto &t_xX\cap  t_yX
\end{diagram}
\leqno (*)
$$
is dominant, where $t_xX,\ t_yX\subseteq\bbbp^N$ are the projective tangent spaces. 

In this paper we show how one can further generalize Ciliberto's result to any
$n$-dimensional subvariety $X$ of $\bbbp^{2n}$. More precisely we show that {\em if\/
$\Tan X=\bbbp^{2n}$ then the map $(*)$ is again surjective,} where $\Tan X$ denotes
the tangent variety of $X$, i.e.\ the closure of the union of all projective tangent
spaces
$t_xX$ at {\em smooth}\/  points of $X$. Note that if $X$ is smooth then by a result
of Fulton and Lazarsfeld \cite{FL} $\Tan X=\bbbp^{2n}$ if and only if $\Sec
X=\bbbp^{2n}$, where $\Sec X$ is the secant variety of
$X$ (see also \cite[1.3 and 4.3.12]{FOV}). As a corollary we obtain
that the St\"uckrad-Vogel cycle
$v_o(X,X)$ has maximal transcendence degree under the assumptions above.   

Throughout this paper we work over the field $K=\bbbc$ of complex numbers. We note
that the results remains valid over any algebraically closed field of characteristic
0 by standard arguments. In the following we will use the notations of our previous
paper \cite{FM2}. As a general reference on the St\"uckrad-Vogel cycle we refer the
reader to \cite{FOV}.

\section{Transcendence degree of the St\"uckrad-Vogel cycle}
Let X be an $n$-dimensional subvariety of $\bbbp^N$ with $N=2n$. Let $x\in X$ be a
generic point and let $U\cong \bbbc^N$ be an affine neighbourhood of $x$ in $\bbbp^N$
such that $x$ corresponds to $0\in \bbbc^N$. After a linear change of coordinates we
can write $X$ in a complex neighbourhood $V$ of $x\hat =0$ as a graph
$$
W\ni u\mapsto (u,f(u))\in V
$$
with $f_u(0)=f(0)=0$, where $W$ is a neighbourhood of $0$ in $\bbbc^N$. 

With these
notation we will show the following explicit criterion as to when
$\Tan X$ has maximal dimension.

\begin{lemma}\label{2.1}
The following are equivalent.

(a) $\Tan X=\bbbp^N$;

(b) $f_{uu}(0)\cdot u$ is a nondegenerate matrix for a generic point $u\in \bbbc^n$.

\end{lemma}

\begin{proof} Assume (a) is satisfied. Consider the map
\begin{diagram}[height=17pt]
W\times\bbbc^n& \cong & T(X\cap V) & \rTo^{can} & \bbbc^N\\
& & (u,\xi) & \rMapsto & (u+\xi, f(u)+f_u(u)\cdot \xi),
\end{diagram}
where $T(X\cap V)$ denotes the tangent bundle of $X\cap V$.
By (a) the differential of this map has maximal rank, i.e.\ the matrix
$$
\left(\begin{matrix}
E_n& E_n\\
f_u(0)+f_{uu}(0)\cdot \xi& f_u(0)
\end{matrix}\right)
=
\left(\begin{matrix}
E_n& E_n\\
f_{uu}(0)\cdot \xi& 0
\end{matrix}\right)
$$
has maximal rank, where $E_n$ is the $n$-dimensional identity matrix.
Hence (b) follows, and by reversing the argument it also follows that (b)
implies (a). 
\end{proof}

We can now formulate the main result of this paper.

\begin{theorem}\label{2.2}
Let $X$ be an $n$-dimensional subvariety of $\bbbp^N$ with $N=2n$. If the tangent
variety\/ $\Tan X$ is equal to $\bbbp^N$, then a generic point $P\in\bbbp^N$ 
is uniquely determined by the ramification locus $R(\pi_P)$ of the linear projection
from
$P$.
\end{theorem}

\begin{proof}
By assumption, $\Tan X=\bbbp^N$ and so in particular $\Sec X=\bbbp^N$. Thus, using
Terracini's lemma \cite[Proposition 4.3.2]{FOV}, for general points
$x,\ y\in X$ the intersection
$t_xX\cap t_yX$ consists of just one point $P=P(x,y)$.  We need to show that the
rational map
$$
\begin{diagram}[height=15pt]
X\times X & \rDashto &  \bbbp^N\\
(x,y) & \rMapsto & P(x,y)
\end{diagram}\leqno (*)
$$
is dominant.  Let $x$ be a general point
of $X$. Clearly it suffices to prove that the rational map 
$$
\begin{diagram}[height=15pt]
X & \rDashto &  T_xX\\
y & \rMapsto & P(x,y)
\end{diagram}\leqno (**)
$$
is dominant. We may assume that $x=[1:0:\ldots :0]$ so that the affine
open set $U=\{x_0=1\}\cong \bbbc^N$ is an affine open neighbourhood of $x$;
note that $x$ then corresponds to the origin in $\bbbc^N$. After a linear change of
coordinates we can write
$X$ in a (complex) neighbourhood of $x\,\hat =\,0$ as a graph 
$$
u\mapsto (u,f(u))\in \bbbc^N
$$
with $f(0)=0$, $f_u(0)=0$. With $y=(u,f(u))$ the point
$$
P=P(0,y)=(p(u),0)\in T_xX,
$$
is the intersection of the linear subspaces $\bbbc^N\times 0$ and
the image of $\xi \mapsto(u+\xi, f(u)+f_u(u)\cdot \xi)$. This leads to
$u+\xi=p(u)$, i.e.\ $ \xi=p(u)-u$ and 
$$
0=f(u)+f_u(u)\cdot\xi=f(u)+f_u(u)(p(u)-u).
$$
By \ref{2.1} the matrix $f_{uu}(0)\cdot u$ is nondegenerate. As 
$$
f_u(u)=f_{uu}(0)\cdot u+\mbox{higher order terms}
$$
the matrix $f_u(u)$ is invertible for a general and sufficiently small point $u$.
Thus we can write
$$
p(u) =u-f_u(u)^{-1} f(u).
$$
We need to show that $u\mapsto p(u)$ has maximal rank. The differential of this map is
$$\begin{array}{rcl}
\eta &\mapsto &
\eta - f_u(u)^{-1} f_u(u)\eta + f_u(u)^{-1}f_{uu}(u)f_u(u)^{-1}\eta\\[2pt]
&&=
f_u(u)^{-1}f_{uu}(u)f_u(u)^{-1}\eta.
\end{array}$$
Thus this map has maximal rank, since $\tilde\eta=f_u(u)^{-1}\eta$ is a 
vector in general position, so by the lemma $f_{uu}(0)\cdot\tilde\eta$ has maximal
rank. As 
$$
f_{uu}(u)\tilde\eta=f_{uu}(0)\tilde\eta +\mbox{higher order terms}
$$
it follows that also $f_{uu}(u)\tilde \eta$ has maximal rank for $u$ sufficiently
small.
\end{proof}

Using \cite[Lemma 3.3]{FM2} we have the following corollary.

\begin{corollary}\label{2.3}
Let $X$ be an $n$-dimensional subvariety of $\bbbp^N$ with $N=2n$. If the tangent
variety\/ $\Tan X$ is equal to $\bbbp^N$, then 
the St\"uckrad Vogel cycle
$v_0=v_0(X,X)$ has maximal transcendence degree.
\end{corollary}

In the case of surfaces in $\bbbp^4$ we recover the result of Ciliberto mentioned in
the introduction.

\begin{corollary} {\em (C.Ciliberto)}
If $X\subseteq \bbbp^4$ is a nondegenerate smooth surface, then a general point
$P\in\bbbp^N$  is uniquely determined by the ramification locus $R(\pi_P)$ of the
linear projection from $P$.
\end{corollary}

\begin{proof}
If $\Tan X\not=\bbbp^4$ then $\dim \Tan X=\dim \Sec X =\dim X+1$ and so by
\cite[4.6.6(1)]{FOV} $\Sec X$ is a linear subspace of dimension 3, which contradicts
the fact that $X$ is nondegenerate. Thus \ref{2.2} implies the result.
\end{proof} 

\begin{remark}
We note that Ciliberto's argument is different from our proof. He
uses the fact that for a nondegenerate surface which is not a cone the second
fundamental form is nondegenerate. 
\end{remark}

\begin{remark}
The example 4.8 given in \cite{FM2} is not a counterexample to the maximality of the
transcendence degree of the St\"uckrad-Vogel cycle as was claimed there. The mistake
comes from the fact that the dimension of the family of lines $\bbbp^1\hto
\bbbp^1\times\bbbp^m$ of bidegree (1,1) was not calculated correctly. It is $2m+1$
(and not $2m-2$), since such a map is given by $(\alpha, j)$, where $\alpha$ is an
automorphism of $\bbbp^1$ and $j:\bbbp^1\hto\bbbp^m$ is the embedding of a line.
\end{remark}

\end{document}